\documentclass{article}
\usepackage[utf8]{inputenc}
\usepackage[left=1in,top=1in,right=1in,bottom=1in]{geometry}
\usepackage[linktocpage=true]{hyperref}
\usepackage{setspace}
\usepackage{amssymb, amsmath, amsthm, graphicx,mathrsfs}
\usepackage{caption,cite}
\usepackage{mathtools,color}

\usepackage{cleveref}
\usepackage{tikz}
\usetikzlibrary{shapes.geometric, backgrounds, positioning, arrows.meta,calc}

\newtheorem{thm}{Theorem}[section]
\newtheorem{cor}[thm]{Corollary}

\newtheorem{defn}{Definition}

\newtheorem{lem}[thm]{Lemma}
\newtheorem{conj}[thm]{Conjecture}
\newtheorem{prop}[thm]{Proposition}
\newtheorem{claim}[thm]{Claim}

\newtheorem{prob}[thm]{Problem}

\renewcommand{\l}{\left}
\renewcommand{\r}{\right}

\newcommand{\sub}{\subseteq}

\title{Generalized Erd\H{o}s-Rogers problems for hypergraphs}
\author{Xiaoyu He \thanks{School of Mathematics, Georgia Institute of Technology, Atlanta, GA 30332. Email: xhe399@gatech.edu.}
\and Jiaxi Nie\thanks{School of Mathematics, Georgia Institute of Technology, Atlanta, GA 30332. Email: jnie47@gatech.edu.}}

\begin{document}

\maketitle

\begin{abstract}
Given $r$-uniform hypergraphs $G$ and $F$ and an integer $n$, let $f_{F,G}(n)$ be the maximum $m$ such that every $n$-vertex $G$-free $r$-graph has an $F$-free induced subgraph on $m$ vertices. We show that $f_{F,G}(n)$ is polynomial in $n$ when $G$ is a subgraph of an iterated blowup of $F$. As a partial converse, we show that if $G$ is not a subgraph of an $F$-iterated blowup and is $2$-tightly connected, then $f_{F,G}(n)$ is at most polylogarithmic in $n$. Our bounds generalize previous results of Dudek and Mubayi for the case when $F$ and $G$ are complete.
\end{abstract}

\section{Introduction}

Given $r$-uniform hypergraphs (henceforth $r$-graphs) $G$ and $F$ and an integer $n \ge 1$, we let $f_{F,G}(n)$ be the maximum integer $m$ such that every $n$-vertex $G$-free $r$-graph contains an $F$-free induced subgraph on $m$ vertices. When $F=K^r_r$ the single edge $r$-graph, determining $f_{K^r_r, G}(n)$ is equivalent to determining the classical off-diagonal hypergraph Ramsey number, which is one of the central problems in extremal combinatorics. Even for graphs, our knowledge of these numbers so far is quite limited: for $K_3$, Ajtai-Komlós-Szemerédi~\cite{ajtai1980note} and Kim~\cite{kim1995ramsey} showed that $f_{K_2,K_3}(n)=\Theta(\sqrt{n\log n})$; for $K_4$, Mattheus and Verstraëte~\cite{mattheus2024asymptotics} showed that $f_{K_2,K_4}(n)=n^{1/3+o(1)}$. We still don't know the correct exponent of $f_{K_2,G}(n)$ when $G$ is $C_4$ or $K_5$.  

Erd\H{o}s and Rogers~\cite{erdos1962construction}, generalizing the off-diagonal Ramsey problem, initiated the study of $f_{K_s,K_t}(n)$; these problems have since attracted significant attention and are known as Erd\H{o}s-Rogers problems. The state of the art on $t = s+1$ are results of Dudek-Mubayi~\cite{dudek2014generalized} and Mubayi-Verstraëte~\cite{mubayi2025order}, establishing the bounds $$\Omega(\sqrt{n\log n/\log\log n})= f_{K_s,K_{s+1}}(n)= O(\sqrt{n}\log n).$$ For $t = s+2$, Sudakov~\cite{sudakov2005large} and Janzer-Sudakov~\cite{janzer2025improved} showed that 
$$
n^{\frac{1}{2}-\frac{1}{6s-6}}(\log n)^{\Omega(1)}= f_{K_s,K_{s+2}}(n)= O(n^{\frac{1}{2}-\frac{1}{8s-10}}(\log n)^3).
$$

Recently, Mubayi-Verstraëte~\cite{mubayi2024erd} and Balogh-Chen-Luo~\cite{balogh2025maximum}, followed soon after by Gishboliner, Janzer and Sudakov~\cite{gishboliner2024induced}, started the systematic study of the function $f_{F,G}(n)$ where $F$ and $G$ are graphs. Their results mostly concern the case when $G$ is a clique, and established bounds for $f_{F,K_r}$ when $F$ satisfies certain properties such as clique-free, bipartite, containing a cycle, or having large minimum degree.

In this paper, we consider the natural generalization of this line of research to hypergraphs. Note that the previous bounds for $f_{F,G}(n)$ are all polynomial. Our first result shows that this is not a coincidence: we find a sufficient condition for $f_{F,G}(n)$ being polynomial, which is satisfied by all pairs of graphs. Let $H$ and $G$ be $r$-graphs. For a vertex $v$ of $H$ and a positive integer $t$, we let $H(v,t)$ be the $r$-graph obtained by adding $t-1$ copies of $v$ to $H$. Further, we let $H(v,F)$ obtained by adding $v(F)-1$ copies of $v$ to $H$, which together with $v$ induce a copy of $F$. 

\begin{defn}
Let $G$ and $F$ be $r$-graphs. We say $G$ is an $F$-iterated blowup if 
\begin{itemize}
    \item[(1)] $G=F$; or
    \item [(2)] $G=H(v,F)$ where $H$ is an $F$-iterated blowup.
\end{itemize}
\end{defn}

\begin{figure}
    \centering
    \includegraphics{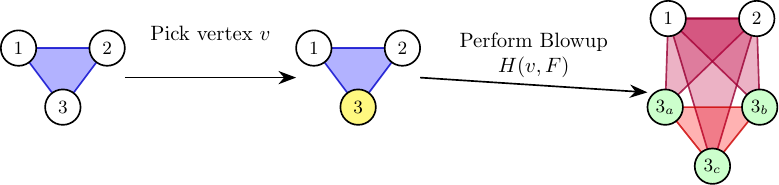}
    \caption{The iterated blowup $H(v,F)$ when $H = F=K_3^3$.}
    \label{fig:iterated blowup}
\end{figure}

When $G=K^3_3$ and $F$ is a subgraph of a $K^3_3$-iterated blowup, Erd\H{o}s and Hajnal~\cite{erdos1972ramsey} showed that $f_{K^3_3,F}(n)\ge n^c$ for some constant $c$ (see also~\cite{conlon2010hypergraph,fox2021independent}). We extend this result, showing that if $G$ is a subgraph of an $F$-iterated blowup, then $f_{F,G}(n)$ has a polynomial lower bound.
\begin{thm}[Proof is in \Cref{section:lower bound}]\label{Theorem:poly lower bounds}
Let $r\ge 2$ and let $G$ and $F$ be $r$-graphs. If $G$ is a subgraph of an $F$-iterated blowup, then there exists a constant $c>0$ depending only on $F, G$ such that, for large enough $n$,
$$
f_{F,G}(n)\ge n^c.
$$
\end{thm}

\Cref{Theorem:poly lower bounds} is a straightforward generalization of a supersaturation argument that was known to Erd\H{o}s (see \cite{fox2021independent}), but we include a proof for completeness.

Note that in the case of graphs, starting from any nonempty graph $F$, one can obtain $F$-iterated blowups with arbitrarily large clique number. Thus, every graph $G$ is a subgraph of an $F$-iterated blowup. Hence \Cref{Theorem:poly lower bounds} reproduces the fact that $f_{F,G}(n)$ is polynomial for graphs. More interesting phenomena appear in hypergraphs, as seen in recent work of Conlon-Fox-Gunby-He-Mubayi-Suk-Verstraëte-Yu~\cite{conlon2024when}, where they prove that if $G$ is tightly connected and not tripartite, then
\begin{equation}\label{equation:Ramsey}
f_{K_3,G}(n)=O((\log n)^{3/2}).    
\end{equation}

Our second result generalizes (\ref{equation:Ramsey}) to the Erd\H os-Rogers setting. We say an $r$-graph $F$ is {\em $k$-tightly connected} if its edges can be ordered as $e_1,~e_2\dots,~e_t$ such that for each $2\le i\le t$ there exists $1\le j\le i-1$ such that $|e_i\cap e_j|\ge k$. In particular, $(r-1)$-tight connectivity is the usual notion of tightly connectivity. For any set $X$, we use $\binom{X}{k}$ to denote the family of all subsets of $X$ of size $k$. The $k$-shadow of an $r$-graph $H$, denoted $\partial_kH:=\bigcup_{e\in E(H)}\binom{e}{k}$, is the $k$-graph whose edges are $k$-subsets of edges of $H$. We use $e(H)$ and $v(H)$ to denote the numbers of edges and vertices in $H$ respectively.

\begin{thm}[Proof is in \Cref{section:not homomorphic}]\label{theorem:not homomorphic}
Let $r\ge 3$ and let $G$ and $F$ be $r$-graphs such that $G$ is 2-tightly connected and is not homomorphic to $F$. Then there exists a constant $c$ depending only on $F$ such that, for large enough $n$,
$$
f_{F,G}(n)\le c(\log n)^{\alpha_F},
$$
where
$$
\alpha_F=\max_{\emptyset\not=F'\sub \partial_2 F}\l\{\frac{e(F')+1}{v(F')-1}\r\}.
$$
\end{thm}

It would be very interesting to characterize pairs $F$ and $G$ such that $f_{F,G}(n)$ is polynomial. In~\cite{conlon2024when}, the following conjecture is proposed.
\begin{conj}[Conjecture 1.1, \cite{conlon2024when}]\label{conjecture:Ramsey}
For a 3-graph $G$, there exists a constant $c=c(G)$ such that $f_{K^3_3,G}(n)\ge n^c$ if and only if $G$ is a subgraph of a $K_3^3$-iterated blowup.
\end{conj}

Extending \Cref{conjecture:Ramsey} , we propose the following.
\begin{conj}\label{conjecture:General}
For any $r$-graphs $F$ and $G$, there exists a constant $c=c(F,G)$ such that $f_{F,G}(n)\ge n^c$ if and only if $G$ is a subgraph of an $F$-iterated blowup.
\end{conj}

Indeed \Cref{theorem:not homomorphic} confirms \Cref{conjecture:General} when $G$ is 2-tightly connected, since it is easy to check that if $G$ is 2-tightly connected and is not homomorphic to $F$, then $G$ is not a subgraph of any $F$-iterated blowup.

Note that when $F=K^3_3$, \Cref{theorem:not homomorphic} only gives $f_{K^3_3,G}(n)\le c(\log n)^2$, which is worse than (\ref{equation:Ramsey}). This is because our proof of \Cref{theorem:not homomorphic} is essentially different from that of (\ref{equation:Ramsey}), in that we sacrifice the exponent to handle a more general class of $F$ and $G$. In particular, when $G$ and $F$ are cliques, say $G=K^r_s$ and $F=K^r_{s+1}$ where $s\ge r\ge 3$, \Cref{theorem:not homomorphic} implies that 
$$f_{K^r_s,K^r_{s+1}}(n)\le c(\log n)^{\frac{\binom{s}{2}+1}{s-1}}.$$ This is much worse than the current best upper bounds of Dudek and Mubayi~\cite{dudek2014generalized}, who show that
\begin{equation}\label{equation:clique}
f_{K^r_s,K^r_{s+1}}(n)\le c(\log n)^{\frac{1}{r-2}}.   
\end{equation}

The method employed by Dudek and Mubayi for (\ref{equation:clique}) is ad-hoc. We provide a new proof of (\ref{equation:clique}), as a consequence of a general upper bound for $G$ and $F$ assuming $G$ is, roughly speaking, far from homomorphic to $F$.

\begin{figure}[h]
    \centering
    \includegraphics{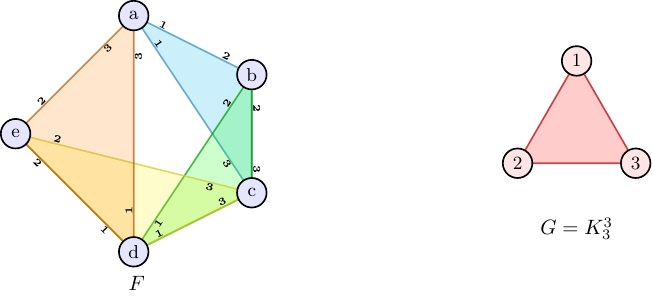}
    \caption{A $2$-shadow homomorphism that is not a homomorphism.}
    \label{fig:shadow-homomorphic}
\end{figure}

\begin{defn}
Given $r$-graphs $F$ and $G$. We say \emph{$F$ is $k$-shadow-homomorphic to} $G$ if we can define for each $S\in \partial_k F$ an $f(S)\in \partial_k G$ and a bijection $g_S:S\rightarrow f(S)$ such that for every edge $e\in E(F)$ there exists an edge $e'\in E(G)$ and a bijection $g: e\rightarrow e'$ such that, for each $S\in \binom{e}{k}$, $g|_{S}=g_S$.
\end{defn}

In other words, $F$ is $k$-shadow-homomorphic to $G$ if we can define bijections from $\partial_k F$ to $\partial_k G$ in a way that glues together consistently along edges of $F$. Note that $1$-shadow-homomorphisms are just homomorphisms; it is also not hard to check that for any $k_1>k_2$, if $F$ is $k_2$-shadow-homomorphic to $G$, then $F$ is $k_1$-shadow-homomorphic to $G$, so in general shadow-homomorphisms are a more permissive notion. For example, let $F$ be the 3-graph with edges $abc,bcd,cde$ and $dea$, then $F$ is $2$-shadow-homomorphic to $K^3_3$ (see \Cref{fig:shadow-homomorphic} for an illustration) but not homomorphic to $K^3_3$.  We remark that shadow homomorphisms are closely related to, but distinct from, the notion of ``pair homomorphism" defined in~\cite{conlon2024off}.

Our last result improves the exponent in \Cref{theorem:not homomorphic} under a more restrictive assumption using shadow homomorphisms.

\begin{thm}[Proof is in \Cref{section:shadow homomorphic}]\label{theorem:shadow-homomorphic}
For $r>k\ge2$, given $r$-graphs $G$ and $F$ such that $G$ is not $k$-shadow-homomorphic to $F$, there exists a constant $c$ depending only on $F$ such that, for large enough $n$,
$$
f^r_{F,G}(n)\le c(\log n)^{\frac{1}{k-1}}.
$$
\end{thm}

For any $s\ge r\ge3$, $K^r_{s+1}$ is not $(r-1)$-shadow-homomorphic to $K^r_s$(Proof is in \Cref{section:shadow homomorphic}), so (\ref{equation:clique}) follows from \Cref{theorem:shadow-homomorphic}. Let $H_t^r$ be the unique $r$-graph with $r+1$ vertices and $t$ edges. For $2\le t\le r$, one can show that $H^r_{t+1}$ is not $(r-1)$-shadow homomorphic to $H^r_{t}$(Proof is in \Cref{section:shadow homomorphic}). Thus, we have the following corollary of \Cref{theorem:shadow-homomorphic}.

\begin{cor}\label{cor:Simplex}
For $r\ge 3$ and $2\le t\le r$, there exists a constant $c$ such that, for large enough $n$,
$$
f^r_{H^r_{t},H^r_{t+1}}(n)\le c(\log n)^{\frac{1}{r-2}}.
$$
\end{cor}

For clarity, we systematically omit all floor and ceiling functions where they are not essential.


\section{Proof of \Cref{Theorem:poly lower bounds}}\label{section:lower bound}

In this section we give the proof of \Cref{Theorem:poly lower bounds}, which is a straightforward generalization of a folklore supersaturation argument that goes back to Erd\H{o}s (see e.g. \cite[Section 6]{fox2021independent}).

\begin{proof}[Proof of \Cref{Theorem:poly lower bounds}]
It suffices to show this assuming $G$ is an $F$-iterated blowup. We will prove this by induction on $G$. 
When $G=F$, the theorem is trivially true. When $G\not=F$, by definition there exists an $F$-iterated blowup $G'$ smaller than $G$ and a vertex $v\in V(G')$ such that $G=G'(v,F)$. By induction, there exists a constant $c'$ such that, for all large enough $n$, $f_{F,G'}(n)\ge n^{c'}$. Let $c_1$ be a constant sufficiently small in terms of $c'$ and $G$. Let $H$ be an $r$-graph such that every $n^{c_1}$-vertex set in $H$ contains a copy of $F$. It suffices to show that $H$ contains a copy of $G$. 

Since $f_{F,G'}(n^{c_1/c'})\ge n^{c_1}$, it follows that every $n^{c_1/c'}$-vertex set in $H$ contains a copy of $G'$. By double counting, the number of copies of $G'$ in $H$ is at least
$$
\frac{\binom{n}{n^{c_1/c'}}}{\binom{n-v(G')}{n^{c_1/c'}-v(G')}}\ge n^{(1-c_1/c')v(G')}.
$$
Note that the number of copies of $G'\setminus v$ in $H$ is at most $n^{v(G')-1}$. Thus there exists a copy of $G'\setminus v$ that can be extended to at least 
$$
n^{(1-c_1/c')v(G')}/n^{v(G')-1}=n^{1-c_1v(G')/c'} \ge n^{c_1}
$$
copies of $G'$ in $H$, as long as $c_1$ is sufficiently small. These extensions together form a copy of $G'(v,n^{c_1})$. By the definition of $H$, the $n^{c_1}$ vertices forming copies of $v$ in $G'(v,n^{c_1})$ contain a copy of $F$. The vertices in this copy of $F$, together with the vertices in the copy of $G'\setminus v$, form a copy of $G'(v,F)=G$, completing the proof.
\end{proof}

\section{Proof of \Cref{theorem:not homomorphic}}\label{section:not homomorphic}
We use the following standard upper tail bound for containing a subgraph in a random graph.
\begin{thm}[Theorem 3.9 from\cite{janson2011random}]\label{theorem:probability of containment}
Let $G(n,p)$ be the Erd\H{o}s-R\'enyi random graph with $n$ vertices and edge probability $p$. Let $F$ be a graph with at least one edge. Then for every sequence of $p=p(n)<1$,
$$
\Pr[F\not\subset G(n,p)]\le \exp(-\Theta(n^{v(F')}p^{e(F')}))
$$
where $F'$ is a non-empty subgraph of $F$ with minimum $n^{v(F')}p^{e(F')}$. 
\end{thm}

\begin{proof}[Proof of \Cref{theorem:not homomorphic}]
Let constants $c_3$ and $c_4$ be sufficiently large; $c_5$ be sufficiently large in terms of $F$; $c_1$ be sufficiently large in terms of $F$ and $c_5$; $c_2$ be sufficiently large in terms of $F$, $c_1$, $c_3$ and $c_4$. We write $\alpha=\alpha_F$ for short.
Let $\ell=c_1\log n$ and let $w=c_2(\log n)^{\alpha}$. 

Consider a random function (coloring) $\beta: \binom{[n]}{2}\rightarrow [\ell]$ where each pair in $\binom{[n]}{2}$ is assigned a color in $[\ell]$ independently and uniformly at random. For each $t\in [\ell]$ we let $G_t$ be the graph on $[n]$ whose edges are all pairs with color $t$ from $\beta$. 

\begin{claim}\label{claim:partialF containment}
With positive probability, for every $t\in [\ell]$ and every $W\sub [n]$ with $|W|\ge w/2$, $G_t[W]$ contains a copy of $\partial_2 F$.     
\end{claim}
\begin{proof}
For every $t\in [\ell]$ and $W\sub [n]$ with $|W|= w/2 $, we let $X_{t,W}$ be the event that $G_t[W]$ is $\partial_2 F$-free. Note that each $G_t$ is a random graph $G(n,p)$ where $p=1/\ell$. Thus by \Cref{theorem:probability of containment}, we have
$$
\Pr[X_{t,W}]\le \exp\l(-\frac{1}{c_3}w^{v(F')}\ell^{-e(F')}\r)
$$
where $F'$ is some non-empty subgraph of $\partial_2F$.

Thus by the union bound, 
$$
\begin{aligned}
\Pr\l[\bigcup_{t\in[\ell], W\in \binom{[n]}{ w/2}} X_{t,W}\r]&\le \ell\binom{n}{w/2}\exp\left(-\frac{1}{c_3}w^{v(F')}\ell^{-e(F')}\right)\\
&\le \exp\l(c_4c_2(\log n)^{1+\alpha}-\frac{c_2^{v(F')}}{c_3c_1^{e(F')}}(\log n)^{v(F')\alpha-e(F')}\r).    
\end{aligned}
$$   
Comparing the exponents of $\log n$ in the two terms above, we note that $v(F')\alpha-e(F')\ge 1+\alpha$ is equivalent to 
$\alpha\ge \frac{e(F')+1}{v(F')-1}$, which is true by definition of $\alpha$. This completes the proof.
\end{proof}

We may thus fix $\beta$ satisfying the conclusion of \Cref{claim:partialF containment}. For each $t\in [\ell]$,take a function $\gamma_t:[n]\rightarrow V(F)$ uniformly at random. We define an $r$-graph $H$ on $[n]$ whose edges are all $r$-tuples $X$ such that there exists $t\in [\ell]$ so that all pairs in $\binom{X}{2}$ are mapped to $t$ by $\beta$ and that $\gamma_t(X)$ is an edge in $F$.

It is not hard to check that $H$ is $G$-free; indeed, since $G$ is 2-tightly connected, a copy of $G$ must have its 2-shadow mapped to the same color $t\in [\ell]$ by $\beta$, and hence $\gamma_t$ will map every edge in this copy of $G$ to an edge in $F$, producing a homomorphism from $G$ to $F$. 

For each $W\sub [n]$ with $|W|=w$, we let $Y_W$ be the event that $W$ is $F$-free. By \Cref{claim:partialF containment}, we know that, for each $W$ with $|W|=w$ and each $t\in [\ell]$, $G_t[W]$ contains at least $\frac{w}{2v(F)}$ vertex-disjoint copies of $\partial_2 F$. The probability that a copy of $\partial_2 F$ in $G_t$ produces a copy of $F$ in $H$ is at least $\frac{v(F)\,!}{v(F)^{v(F)}}$. Thus
$$
\Pr[Y_W]\le \l(1-\frac{v(F)\,!}{v(F)^{v(F)}}\r)^{\frac{\ell w}{2v(F)}}.
$$
By the union bound,
$$
\Pr\l[\bigcup_{W\in\binom{[n]}{w}} Y_W\r]\le \binom{n}{w}\l(1-\frac{v(F)\,!}{v(F)^{v(F)}}\r)^{\frac{\ell w}{2v(F)}}\le \exp\l(\log n\cdot w-\frac{1}{c_5}\ell w\r)=\exp\l((1-\frac{c_1}{c_5})\log n\cdot w\r)<1.
$$
This means that, with positive probability, for each $W\sub [n]$ with $|W|=w$, $H[W]$ contains a copy of $F$.
\end{proof}

\section{Proof and applications of \Cref{theorem:shadow-homomorphic}}\label{section:shadow homomorphic}

In this section, we prove \Cref{theorem:shadow-homomorphic} which improves the bound for $f^r_{F,G}$ under the more restrictive condition that $G$ is not $k$-shadow-homomorphic to $F$. 

We make use of the extended form of the Janson Inequality, in the following form.

\begin{thm}[Theorem 8.1.2, \cite{alon2016probabilistic}]\label{theorem:Extended Janson}
Let $\Omega$ be a finite set, and let $R$ be a random subset of $\Omega$ given by $\Pr[r\in R]=p_r$, these events being mutually independent over $r\in\Omega$. Let $\{A_i\}_{i\in I}$ be a finite collection of subsets of $\Omega$. Let $B_i$ be the event $A_i\subseteq R$. For $i,j\in I$, we write $i\sim j$ if $i\not=j$ and $A_i\cap A_j\not=\emptyset$. Define
$$
\Delta=\sum_{i\sim j}\Pr[B_i\cap B_j],
$$
where the sum is over ordered pairs, and
$$
\mu=\sum_{i\in I}\Pr[B_i].
$$
If $\Delta\ge \mu$, then
$$
\Pr[\cap_{i\in I}\overline{B_i}]\le \exp(-\mu^2/\Delta).
$$
\end{thm}

\begin{proof}[Proof of \Cref{theorem:shadow-homomorphic}]
Let constants $c_3,c_4$ be sufficiently large in terms of $F$; $c_1$ be sufficiently large in terms of $c_3$ and $c_4$; $c_2$ be sufficiently large in terms of $c_1$.

Let $w=c_2(\log n)^{\frac{1}{k-1}}$. We want to construct an $n$-vertex $G$-free $r$-graph $H$ such that every $w$-vertex set in $H$ contains a copy of $F$. We construct $H$ randomly as follows. Let $[n]$ be the vertex set of $H$. For each $S\in \binom{[n]}{k}$, we take an $f(S)\in\partial_k F$ uniformly at random and then take a bijection $g_S:S\rightarrow f(S)$ uniformly at random. For any $X\in \binom{[n]}{r}$, we let $X$ be an edge of $H$ if and only if there exist an edge $e\in E(F)$ and a bijection $g:X\rightarrow e$ such that, for each $S\in \binom{X}{k}$, $g|_S=g_S$. It is not hard to check that $H$ is $G$-free because otherwise the functions $g_S$ would glue together to a $k$-shadow-homomorphism from $G$ to $F$, which is a contradiction. 
\begin{claim}\label{claim:F-free}
Let $W\in \binom{[n]}{w}$. The probability that $H[W]$ is $F$-free is at most
$$
\exp\l(-\frac{1}{c_1}w^k\r).
$$
\end{claim}
\begin{proof}
To simplify the analysis we will partition $H[W]$ into $v(F)$ parts and only consider transversal copies of $F$.

Let $v_1,v_2,\dots,v_t$ be all vertices in $F$. Consider a partition $W=V_1\sqcup V_2\sqcup\dots\sqcup V_t$ where $|V_i|=w/t$ for each $i$. We say a $k$-set $S=\{s_1,s_2,\dots,s_k\}\subseteq W$ is \emph{good} if there exists $\{v_{i_1},v_{i_2},\dots,v_{i_k}\}\in \partial_k F$ such that $s_t\in V_{i_t}$ for all $1\le t\le k$. For each good $k$-set $S=\{s_1,s_2,\dots,s_k\}\subseteq W$, we say $S$ is \emph{faithful} if $g_S(s_t)=v_{i_t}$ for all $1\le i\le k$. Let $Y_S$ be the event that $S$ is faithful; note that these events are mutually independent over all good $k$-sets $S\subseteq W$.

We say a set $X=\{x_1,x_2,\dots,x_t\}\subseteq W$ is transversal if $x_i\in V_i$ for each $1\le i\le t$. For each transversal $t$-set $X=\{x_1,x_2,\dots,x_t\}$, we let $Y_X$ be the event that, for any good $k$-set $S\subseteq X$ and any $1\le i\le t$, $g_S(x_i)=v_i$ if $x_i\in S$; in other words, all good $k$-sets $S\subseteq X$ are faithful. Note that $Y_X$ happens if and only if $H[X]$ forms a copy of $F$ in the right order. Thus, it suffices to prove an upper bound for the probability that none of the events $Y_X$ happens. 

Let $\Omega$ be the set of all good $k$-sets in $W$ and let $R$ be the random subset of $\Omega$ consisting of all faithful good $k$-sets. Then $Y_X$ can be viewed as the event that all good $k$-sets $S\subseteq X$ are contained in $R$. Let 
$$\mu=\sum_X\Pr[Y_X]$$
where $X$ ranges over all transversal $t$-sets, and let
$$
\Delta=\sum_{Y_{X_1}\sim Y_{X_2}}\Pr[Y_{X_1}\cap Y_{X_2}]
$$
where $Y_{X_1}\sim Y_{X_2}$ means $X_1\cap X_2$ contains a good $k$-set, and in particular, $|X_1\cap X_2|\ge k$. It is not hard to check that $\mu\ge \frac{1}{c_3}w^t$, $\Delta\le c_4w^{2t-k}$ and that $\Delta\ge \mu$ given that $n$ is sufficiently large. By \Cref{theorem:Extended Janson}(the Extended Janson Inequality), we have
\[
\Pr\l[\bigcap_{X}\overline{Y_X}\r]\le \exp\l(-\frac{\mu^2}{2\Delta}\r)\le \exp\l(-\frac{1}{c_1}w^k\r),
\]
as desired.
\end{proof}

By the union bound and \Cref{claim:F-free}, the probability that there exists a $w$-vertex set $W$ in $H$ such that $H[W]$ is $F$-free is at most
$$
\binom{n}{w}\exp\l(-\frac{1}{c_1}w^k\r)\le \exp\l(\log n\cdot w-\frac{1}{c_1}w^k\r)=\exp\l(\l(c_2-\frac{c_2^k}{c_1}\r)(\log n)^{\frac{k}{k-1}}\r)<1.
$$
Thus, with positive probability, $H$ satisfies the desired properties. This completes the proof of \Cref{theorem:shadow-homomorphic}.
\end{proof}

Next, we check the non-shadow-homomorphisms that we claimed in the introduction. We first prove a useful lemma.

\begin{lem}\label{lemma:injectivity}
Let $G$ and $F$ be two $r$-graphs such that $G$ is $(r-1)$-shadow-homomorphic to $F$, so that for each $S\in\partial_{r-1}G$ there is an $f(S)\in\partial_{r-1}F$ and a bijection $g_S:S\rightarrow f(S)$, and for each $E\in E(G)$ there is an $f(E)$ and a bijection $g_E:E\rightarrow f(E)$ such that, for any $S\in E$, $g_S=g_E|S$. Let $E_1,E_2$ and $E_3$ be three edges of $G$ such that $|E_1\cap E_2|=|E_1\cap E_3|=|E_2\cap E_3|=r-1$ and $|E_1\cap E_2\cap E_3|=r-2$. Then $f(E_1)\not=f(E_2)$. 
\end{lem}

\begin{proof}
Let $S=E_1\cap E_2$, $E_1=S\cup\{v_1\}$ and $E_2=S\cup \{v_2\}$. Suppose for contradiction that $f(E_1)=f(E_2)$. Then by definition, for each $v\in S$, $g_{E_1}(v)=g_S(v)=g_{E_2}(v)$, and hence we have $g_{E_1}(v_1)=g_{E_2}(v_2)$. Note that $v_1,v_2\in E_3$. By definition, $g_{E_3}(v_1)=g_{E_3\cap E_1}(v_1)=g_{E_1}(v_1)$, and similarly, $g_{E_3}(v_2)=g_{E_2}(v_2)$. Thus $g_{E_3}(v_1)=g_{E_3}(v_2)$, which contradicts the fact that $g_{E_3}$ is a bijection.
\end{proof}

\begin{prop}
For all $s\ge r\ge 2$, $K^r_{s+1}$ is not $(r-1)$-shadow-homomorphic to $K^r_{s}$.
\end{prop}

\begin{proof}
Let $V=\{v_1,~v_2,\dots,v_s\}$ and $U=\{u_1,~u_2,\dots,u_{s+1}\}$ be the vertex sets of $K^r_{s}$ and $K^r_{s+1}$ respectively.
Suppose for contradiction that $K^r_{s+1}$ is $(r-1)$-shadow-homomorphic to $K^r_{s}$, so that we can pick for each $S\in \binom{U}{r-1}$ an $f(S)\in \binom{V}{r-1}$ and a bijection $g_S:S\rightarrow f(S)$, and for each $E\in\binom{U}{r}$ an $f(E)\in \binom{V}{r}$ and a bijection $g_E:E\rightarrow f(E)$ such that, for any $S\sub E$, $g_S=g_E|_S$. 

Let $U'=\{u_1,~u_2,\dots,u_{r-1}\}$. Without loss of generality, we may assume $g_{U'}(u_i)=v_i$ for every $1\le i\le r-1$. Note that for each $r\le j\le s+1$, $g_{U'\cup\{u_j\}}(u_j)\not\in V':=\{v_1,~v_2,\dots,v_{r-1}\}$; thus $g_{U'\cup\{u_j\}}(u_j)\in V\setminus V'$. By \Cref{lemma:injectivity}, for $r\le j_1<j_2\le s+1$, $f(U'\cup\{u_{j_1}\})\not=f(U'\cup\{u_{j_2}\})$ and hence $g_{U'\cup\{u_{j_1}\}}(u_{j_1})\not=g_{U'\cup\{u_{j_2}\}}(u_{j_2})$. Thus, the set $\{g_{U'\cup\{u_j\}}(u_j):r\le j\le s+1\}$ consists of $s-r+2$ distinct vertices in $V\setminus V'$. But $|V\setminus V'|=s-r+1$, so this is a contradiction.
\end{proof}

Thus, by \Cref{theorem:shadow-homomorphic}, $f_{K^r_s,K^r_{s+1}}\le c(\log n)^{\frac{1}{r-2}}$. This matches (\ref{equation:clique}). 

\begin{prop}\label{proposition:simplex}
For $2\le t\le r$, $H^r_{t+1}$ is not $(r-1)$-shadow-homomorphic to $H^r_{t}$.
\end{prop}

\begin{proof}
Let $V=\{v_1,v_2,\dots,v_{r+1}\}$ and $U=\{u_1,u_2,\dots,u_{r+1}\}$ be the vertex sets of $H^r_{t}$ and $H^r_{t+1}$ respectively. Suppose for contradiction that $H^r_{t+1}$ is $(r-1)$-shadow-homomorphic to $H^r_{t}$, then by definition, we can define for each $S\in \partial_{r-1}H^r_{t+1}$ an $f(S)\in\partial_{r-1}H^r_{t}$ and a bijection $g_S:S\rightarrow f(S)$, and define for each $E\in E(H^r_{t+1})$ an $f(E)\in E(H^r_{t})$ such that for any $S\subset E$, $g_S=g_E|S$.

Since we are working with $r$-uniform hypergraphs on only $r+1$ vertices, every triple of edges of $H_{t+1}^r$ satisfy the conditions of \Cref{lemma:injectivity}. Thus, by \Cref{lemma:injectivity}, the function $f:E(H^r_{t+1})\rightarrow E(H^r_{t})$ is injective. But $|E(H^r_{t+1})|=t+1>t=|E(H^r_{t})|$, so this is a contradiction.
\end{proof}

\Cref{cor:Simplex} follows immediately from \Cref{theorem:shadow-homomorphic} and \Cref{proposition:simplex}.

\section*{Concluding Remarks}
A direct generalization of the proof of (\ref{equation:Ramsey}) in~\cite{conlon2024when} would give the following theorem.
    \begin{thm}
Let $G$ and $F$ be $3$-graphs such that $G$ is tightly connected and is not homomorphic to $F$. If $\frac{e(F')}{v(F')-1}\le \frac{e(\partial_2F)}{v(F)-1}$ for any nonempty $F'\subseteq\partial_2F$, then there exists a constant $c$ depending only on $F$ such that, for large enough $n$,
$$
f_{F,G}(n)\le c(\log n)^{\frac{e(\partial_2F)}{v(F)-1}}.
$$
\end{thm}
This theorem is strictly stronger than \Cref{theorem:not homomorphic} whenever it applies. We believe the extra restriction on $F$ is not necessary. 

\begin{conj}
Let $G$ and $F$ be $3$-graphs such that $G$ is tightly connected and is not homomorphic to $F$. Then there exists a constant $c$ depending only on $F$ such that, for large enough $n$,
$$
f_{F,G}(n)\le c(\log n)^{\beta_F},
$$
where
$$
\beta_F=\max_{\emptyset\not=F'\sub \partial_2 F}\l\{\frac{e(F')}{v(F')-1}\r\}.
$$
\end{conj}

Determining the magnitude of $f^3_{H^3_{2},H^3_{3}}(n)$ seems to be (in some sense) the minimum non-trivial question of this kind.
By \Cref{cor:Simplex} with $r=3$ and $t=2$, we have $f^3_{H^3_{2},H^3_{3}}(n)\le c\log n$. On the other hand, from the Ramsey result of $H^3_{3}$~\cite{fox2021independent}, we know that $f^3_{H^3_{2},H^3_{3}}(n)\ge f^3_{K^3_{3},H^3_{3}}(n)\ge c\frac{\log n}{\log\log n}$. We are not sure which bound is closer to the truth.

\begin{prob}
Determine the magnitude of $f^3_{H^3_{2},H^3_{3}}(n)$.
\end{prob}

Finally, regarding the clique Erdős–Rogers problem, we would like to highlight a problem posed by Conlon, Fox, and Sudakov~\cite{conlon2015recent}.
\begin{prob}
Is it the case that $f_{K^4_s,K^4_{s+1}}=(\log n)^{o(1)}$ for every $s\ge 4$?
\end{prob}
It seems that all the methods in this paper, and previous work on this topic, meet a natural barrier at $(\log n)^c$, so entirely new constructions will be necessary to settle this problem in the affirmative.

\section*{Acknowledgement}
The authors would like to thank Dhruv Mubayi for suggesting the topic of this paper.

\bibliographystyle{abbrv}
\bibliography{Refs}

\end{document}